\documentclass[12pt]{article}
\usepackage[dvips]{graphicx}
\usepackage{amssymb,amsfonts,amsmath}
\usepackage{amsmath,amsmath,amsfonts}
\usepackage{amssymb}

\usepackage{epsfig}

\begin{document}

\begin{center}
{\Large \bf Non-Bayesian particle filters}
$ $\\
$ $\\
{\bf\large Alexandre J.\ Chorin and Xuemin Tu}\\
$ $\\
Department of Mathematics, University of California at Berkeley\\
and\\
Lawrence Berkeley National Laboratory\\
Berkeley, CA, 94720
\end{center}

\begin{abstract}
Particle filters for data assimilation in nonlinear problems use
``particles" (replicas of the underlying system) to generate a sequence of
probability density functions (pdfs) through a Bayesian process. 
This can be expensive because a significant number of particles 
has to be used to maintain accuracy.
We offer here an alternative, in which the relevant pdfs are sampled directly
by an iteration.
An example is discussed in detail. 
\end{abstract}

{\bf Keywords} {particle filter, chainless sampling,  normalization factor,  iteration,  non-Bayesian}

\section{Introduction.}
{T}here are many problems in science in which the state
of a system must be identified from an uncertain equation supplemented by a stream of noisy data (see e.g. \cite{Do2}). A natural model of this situation consists of a stochastic differential equation (SDE):
\begin{equation}
d{\bf x}={\bf f}({\bf x},t) \, dt+ g({\bf x},t) \, d{\bf w},
\label{eq:datass}
\end{equation}
where ${\bf x}=(x_1,x_2,\dots,x_m)$ is an $m$-dimensional vector, $d{\bf w}$ is $m$-dimensional Brownian motion, ${\bf f}$ is an $m$-dimensional vector function, and $g$ is a
scalar (i.e., an $m$ by $m$ diagonal matrix of the form $gI$, where $g$ is a scalar
and $I$ is the identity matrix). The Brownian motion encapsulates all the uncertainty in this
equation. The initial state
${\bf x}(0)$ is assumed given and may be random as well.

As the experiment unfolds, it is observed, and the
values ${\bf b}^n$ of a measurement process are recorded at times $t^n$;
for simplicity assume $t^n=n\delta$, where $\delta$ is a fixed
time interval and $n$ is an integer.
The measurements are related to the evolving state ${\bf x}(t)$ by
\begin{equation}
{\bf b}^n={\bf h}({\bf x}^n)+G{\bf W}^n,
\label{eq:observe}
\end{equation}
where ${\bf h}$ is a $k$-dimensional, generally nonlinear, vector function with $k \le m$, $G$ is a diagonal matrix, ${\bf x}^n={\bf x}(n\delta)$,  and ${\bf W}^n$ is a vector whose
components are independent Gaussian variables of mean 0 and variance 1, independent also of the Brownian motion in equation (\ref{eq:datass}). The task is to estimate ${\bf x}$ on the basis of equation (\ref{eq:datass}) and the observations (\ref{eq:observe}).

If the system~(\ref{eq:datass}) is linear and the data are Gaussian, the solution can be found via the Kalman-Bucy filter. In the general case, it is natural to try to estimate ${\bf x}$ as the mean of its evolving probability density.
The initial state ${\bf x}$ is known and so is its probability density; all one has to do is
evaluate sequentially the density $P_{n+1}$ of ${\bf x}^{n+1}$ given the probability
density $P_n$ of
${\bf x}^n$ and the data ${\bf b}^{n+1}$. This can be done by following
``particles" (replicas of the system) whose empirical distribution approximates $P_n$.
In a Bayesian filter (see e.g 
\cite{Mac1,LS,Dou1,Ar2,Gi1,Ch11,Dow1,toapp},
one uses the pdf $P_n$ and equation (\ref{eq:datass}) to generate
a prior density, and then one uses the new data ${\bf b}^{n+1}$ to generate a
posterior density $P_{n+1}$. In addition, one may have to sample backward to take into account the information
each measurement provides about the past and avoid having too many identical particles. Evolving particles is typically expensive, and the backward sampling,
usually done by Markov chain Monte Carlo (MCMC), can be expensive as well,
because the number of particles needed can grow catastrophically
(see e.g. 
\cite{Sny}). 

In this paper we
offer an alternative to the standard approach, in which $P_{n+1}$ is sampled directly without recourse to Bayes' theorem and backward sampling, if needed, is done
by chainless Monte Carlo \cite{Ch101}. 
Our direct sampling 
is based on a representation of a variable with density $P_{n+1}$ by a collection of functions of  Gaussian
variables parametrized by the support of $P_n$, with parameters found by iteration. 
The construction is related to chainless sampling as described in \cite{Ch101}.
The idea in chainless sampling is to produce a sample of a large set of variables by sequentially sampling
a growing sequence of nested conditionally independent subsets.
As observed in \cite{We1,We2}, chainless sampling for a SDE reduces to interpolatory sampling, as explained below. 
Our construction will be explained in the following sections through an example where the position of a ship is deduced from the
measurements of an azimuth, already used as a test bed in \cite{Gi1,Gordon,Carpenter}.

\section{Sampling by interpolation and iteration.}

First we explain how to sample via interpolation and iteration in a simple
example, related to the example and the construction in \cite{We1}. Consider the scalar SDE
\begin{equation}
dx=f(x,t)dt+ \sqrt{\sigma}\,dw;
\label{scalar}
\end{equation}
we want to find sample paths $x=x(t), 0\le t \le 1$,
subject to the conditions $x(0)=0,x(1)=X$.

Let $N(a,v)$ denote a Gaussian variable with mean $a$ and variance $v$.
We first discretize equation (\ref{scalar}) on a regular mesh
$t^0,t^1,\dots,t^N$,
where 
$t^n=n\delta$, $\delta=1/N$, $0 \le n \le N$,
with $x^n=x(t^n)$, and, following \cite{We1},
use a balanced implicit discretization \cite{MPS,KP}:
$$x^{n+1}=x^n+f(x^n,t^n)\delta+(x^{n+1}-x^n)f'(x^n)\delta+W^{n+1},$$
where $f'(x^n,t^n)=\frac{\partial f}{\partial x^n}(x^n,t^n)$ 
 and $W^{n+1}$ is $N(0,\sigma/N)$.
The joint probability density of the variables $x^1,\dots,x^{N-1}$ is
$Z^{-1}\exp(-\sum_0^N V^i)$, where $Z$
is the normalization constant and 
\begin{align*}
V_i&=\frac{\left((1-\delta f')(x^{n+1}-x^n)-\delta f\right)^2}{2\sigma\delta}\\
&=\frac{\left(x^{n+1}-x^n-\delta f/(1-\delta f')\right)^2}{2\sigma_n},
\end{align*}
where $f,f'$ are functions of the $x^j$, and 
$\sigma_n=\sigma\delta/(1-\delta f')^2$ (see \cite{St10}). One can obtain sample solutions by sampling this density, e.g. by MCMC,
or one can obtain them by interpolation (chainless sampling), as follows.

Consider first the special case $f(x,t)=f(t)$, so that in particular
$f'=0$. Each increment $x^{n+1}-x^n$ is
now a $N(a_n,\sigma/N)$ variable, with the $a_n=f(t^n)\delta$ known explicitly. Let $N$ be a power of $2$.  Consider the variable
$x^{N/2}$. On one hand,
$$x^{N/2}=\sum_1^{N/2}(x^n-x^{n-1})=N(A_1,V_1),$$
where $A_1=\sum_1^{N/2}a_n,
V_1=\sigma/2$.
On the other hand, 
$$X=x^{N/2}+\sum_{N/2+1}^{N}(x^n-x^{n-1}),$$ 
so that
$$x^{N/2}=N(A_2,V_2),$$ 
with $$A_2=X-\sum_{N/2+1}^{N-1}a_n, \quad V_2=V_1.$$ 
The pdf of $x^{N/2}$ is the product of the two pdfs; one can check that
\begin{align*}
&\exp\left(-\frac{(x-A_1)^2}{2V_1}\right)\exp\left(-\frac{(x-A_2)^2}{2V_2}\right)\\
=&\exp\left(-\frac{(x-\bar{a})^2}{2\bar{v}}\right)\exp(-\phi),
\end{align*}
where $\bar{v}=\frac{V_1V_2}{V_1+V_2}$,
$\bar{a}=\frac{V_1A_1+V_2A_2}{V_1+V_2}$, and $\phi=\frac{(A_2-A_1)^2}{2(V_1+V_2)}$;
$e^{-\phi}$ is the probability of getting from the origin to $X$, up to a normalization 
constant.

Pick a sample $\xi_1$ from the
$N(0,1)$ density; one obtains a sample of $x^{N/2}$ by 
setting $x^{N/2}=\bar{a}+\sqrt{\bar{v}}\xi_1$. 
Given a sample of $x^{N/2}$ one can similarly sample $x^{N/4},x^{3N/4}$, then $x^{N/8}$,
$x^{3N/8}$, etc., until all the $x^j$ have been sampled.
If we define ${\bf \xi}=(\xi_1,\xi_2,\dots,\xi_{N-1})$, then for
each choice of ${\bf \xi}$ we find a sample $(x^1,\dots,x^{N-1})$ such that
\begin{align}
&\exp\left(-\frac{\xi_1^2+\cdots+\xi_{N-1}^2}{2}\right)\exp\left(-\frac{(X-\sum_n a_n)^2}{2\sigma}\right)\nonumber\\
=&\exp\left(-\frac{(x^1-x^0-a_0)^2}{2\sigma/N}-\frac{(x^2-x^1-a_1)^2}{2\sigma/N}\right.\nonumber\\
&\left.-\dots -\frac{(x^N-x^{N-1}-a_{N-1})^2}{2\sigma/N}\right),
\label{palim}
\end{align}
where the factor $\exp\left(-\frac{(X-\sum_n a_n)^2}{2\sigma}\right)$ on the left 
is the probability of
the fixed end value $X$ up to a normalization constant. 
In this linear problem, this factor is the same for all the samples and therefore harmless.
One can repeat this sampling process for multiple choices of the variables $\xi_j$;  
each sample of the corresponding set of $x^n$ is independent of any
previous samples of this set.  

Now return to the general case. The functions $f$, $f'$ are now functions of the $x^j$.
We obtain a sample of the probability density
we want by iteration. 
First pick $\Xi=(\xi_1,\xi_2,\dots,\xi_{N-1})$, where
each $\xi_j$ is drawn independently from the $N(0,1)$ 
density (this vector remains fixed during the iteration). 
Make a first guess ${\bf x}^0=(x^1_0,x^2_0,\dots,x^{N-1}_0)$ 
(for example, if $X\ne 0$, pick ${\bf x}=0$).
Evaluate the functions $f,f'$ at ${\bf x}^j$ 
(note that now
$f'\ne 0$, and therefore the variances of the various increments
are no longer constants). We are back in previous case, and can find 
values of the increments $x^{n+1}_{j+1}-x^n_{j+1}$ corresponding to the values 
of $f,f'$ we have. Repeat the process starting with the new iterate. 
If the vectors ${\bf x}^j$ converge to a vector 
${\bf x}=(x^1,\dots,x^{N-1})$, 
we obtain, in the limit, equation 
(\ref{palim}), where now on the right side $\sigma$ depends on $n$ so 
that $\sigma=\sigma_n$, and both $a_n,\sigma_n$ are functions of the 
final ${\bf x}$. The left hand side of (\ref{palim}) becomes:
\[
\exp\left(-\frac{\xi_1^2+\cdots+\xi_{N-1}^2}{2}\right)\exp\left(-\frac{(X-\sum_n a_n)^2}{2\sum_n\sigma_n}\right).
\]
Note that now the factor $\exp\left(-\frac{(X-\sum_n a_n)^2}{2\sum_n\sigma_n}\right)$ is 
different from 
sample to sample, and changes the relative weights of the different samples. 
In averaging, one should take this 
factor as weight, or resample as described at the end of the following section. In order to obtain more uniform 
weights, one also can use the strategies in \cite{Ch101, We1}.

One can readily see that the
iteration converges if $KTM<1$, where $K$ is the Lipshitz constant
of $f$, $T$ is the length of the interval on which one works (here $T=1$),
and $M$ is the maximum norm of the
vectors ${\bf x}^{j+1}-{\bf x}^j$. If this inequality is not satisfied for the iteration above, 
it can be re-established by a suitable underrelaxation. 
One should 
course choose $N$ large enough so that the results are converged in $N$.
We do not provide more details here because they are extraneous to our
purpose, which is to explain chainless/interpolatory sampling and the use of reference
variables in a simple context.

\section{The ship azimuth problem.}

The problem we focus on is discussed in \cite{Gi1,Gordon,Carpenter},  where it is used to demonstrate
the capabilities of particular Bayesian filters. 
A ship sets out from a point $(x_0,y_0)$
in the plane and undergoes a random walk,
\begin{align}
x^{n+1}&=x^n+dx^{n+1},\nonumber\\
y^{n+1}&=y^n+dy^{n+1},
\label{eq1}
\end{align}
for $n \ge 0$, and with $x^0=y^0$
given, and $dx^{n+1}=N(dx^n,\sigma)$, $dy^{n+1}=N(dy^n,\sigma)$,
i.e., each displacement is a sample of a Gaussian random variable whose
variance $\sigma$ does not change from step to step and whose mean is the value of the
previous displacement.
An observer makes noisy measurements of the azimuth $\arctan (y^n/x^n)$, 
recording
\begin{equation}
b^n=\arctan\frac{y^n}{x^n}+N(0,s).
\label{eq2}
\end{equation}
where the variance $s$ is also fixed; here the observed quantity $b$ is scalar and is not be denoted by a boldfaced letter. 
The problem is to reconstruct the positions
${\bf x}^n=(x^n,y^n)$ from equations (\ref{eq1},\ref{eq2}).
We take the same parameters as \cite{Gi1}: $x_0=0.01,y_0=20$,
$dx^1=0.002$, $dy^1=-0.06$, $\sigma=1\cdot10^{-6}, s=25\cdot10^{-6}$.
We follow numerically $M$ particles,
all starting from $X_i^0=x_0,Y_i^0=y_0$,
as described in the following sections, and we
estimate the ship's position at time $n\delta$ as the mean of the locations
${\bf X}^n_i=(X^n_i,Y^n_i),i=1,\dots,M$ of the particles at that time.
The authors of \cite{Gi1} also show numerical results for runs with varying
data and constants;  we discuss those refinements in section 6 below.

\section{Forward step.}
Assume we have a collection of $M$ particles ${\bf X}^n$ at time
$t^n=n\delta$ whose empirical density approximates $P_n$;
now we find increments $d{\bf X}^{n+1}$ such that the empirical
density of ${\bf X}^{n+1}={\bf X}^n+d{\bf X}^{n+1}$ approximates $P_{n+1}$. $P_{n+1}$ is
known implicitly: it is the product of the density that can be deduced from the SDE and the one that comes from the observations, with the appropriate normalization.
If the
increments were known, their probability $p$ 
(the density $P_{n+1}$ evaluated at the resulting positions ${\bf X}^{n+1}$)
would be known,
so $p$ is a function of $d{\bf X}^{n+1}$, ~$p=p(d{\bf X}^{n+1})$. For each particle $i$, 
we are going to sample a Gaussian reference density, obtain a sample of probability $\rho$, then solve (by iteration)
the equation
\begin{equation}
\rho=p(d{\bf X}^{n+1}_i)
\label{main}
\end{equation}
to obtain $d{\bf X}_i^{n+1}$.

Define 
$f(x,y)=\arctan(y/x)$ and $f^n=f(X^n,Y^n)$.
We are working on one particle at a time, so the
index $i$ can be temporarily suppressed.
Pick two independent samples $\xi_x$, $\xi_y$ from
a $N(0,1)$ density 
(the reference density in the present calculation), 
and set $\rho=\frac{1}{2\pi}\exp\left(-\frac{\xi_x^2}{2}-\frac{\xi_y^2}{2}\right)$; 
the variables $\xi_x$, $\xi_y$ remain unchanged until the end of the iteration. We are 
looking for displacements $dX^{n+1}$, $dY^{n+1}$, and  
parameters $a_x,a_y,v_x,v_y, \phi$, such that:
\begin{align}
2\pi\rho=&\exp\left(-\frac{(dX^{n+1}-dX^n)^2}{2\sigma}-\frac{(dY^{n+1}-dY^n)^2}{2\sigma}
\right.\nonumber\\
&\left.-\frac{(f^{n+1}-b^{n+1})^2}{2s}\right)\exp(\phi)\nonumber\\
=&\exp\left(-\frac{(dX^{n+1}-a_x)^2}{2v_x}-\frac{(dY^{n+1}-a_y)^2}{2v_y}\right)
\label{forward} 
\end{align}
The first equality states what we wish to accomplish: find increments $dX^{n+1}$, $dY^{n+1}$,
functions respectively of $\xi_x,\xi_y$, whose probability with respect to $P_{n+1}$
is $\rho$. The factor $e^{\phi}$ 
is needed to
normalize this term ($\phi$ is called below
a ``phase"). The second equality says how the goal is reached: 
we are looking for parameters  $a_x,a_y,v_x,v_y,$ (all functions of ${\bf X}^n$) such that the increments
are samples of Gaussian variables with these parameters,  with the assumed probability. 
One should remember that in our example
the mean of $dX^{n+1}$ is $dX^n$, and similarly for $dY^{n+1}$.
We are not representing $P_{n+1}$ as a function of a single Gaussian- there is a different Gaussian
for every value of ${\bf X}^n$.

To satisfy the second equality we set 
up an iteration for vectors $d{\bf X}^{n+1,j}(=d{\bf X}^j$ for brevity) that converges
to $d{\bf X}^{n+1}$. Start with $d{\bf X}^0=0$. We now explain how to compute
$d{\bf X}^{j+1}$ given $d{\bf X}^j$.

Approximate the observation equation (\ref{eq2}) by
\begin{equation}
f({\bf X}^j)+f_x\cdot(dX^{j+1}-dX^j)+f_y\cdot (dY^{j+1}-dY^j) =b^{n+1}+N(0,s),
\label{obs}
\end{equation}
where the derivatives $f_x,f_y$ are, like $f$,
evaluated at ${\bf X}^j={\bf X}^n+d{\bf X}^j$,
i.e., approximate the observation equation by its Taylor series expansion
around the previous iterate.
Define a variable
$\eta^{j+1}=(f_x\cdot dX^{j+1}+f_y \cdot dY^{j+1})/\sqrt{f_x^2+f_y^2}$.
The approximate observation equation says that $\eta^{j+1}$ is a $N(a_1,v_1)$ variable,
with
\begin{align}
a_1&=-\frac{f-f_x\cdot dX^j-f_y\cdot dY^j-b^{n+1}}{\sqrt{f_x^2+f_y^2}}\, ,\nonumber\\
v_1&=\frac{s}{f_x^2+f_y^2}.
\end{align}
On the other hand, from the equations
of motion one finds that $\eta^{j+1}$ is $N(a_2,v_2)$,
with $a_2=(f_x\cdot dX^n+f_y\cdot dY^n)/\sqrt{f_x^2+f_y^2}$ and $v_2=\sigma$. 
Hence the pdf of $\eta^{j+1}$ is, up to normalization factors, 
$$
\exp\left(-\frac{(x-a_1)^2}{2v_1}-\frac{(x-a_2)^2}{2v_2}\right)
=\exp\left(-\frac{(x-\bar{a})^2}{2\bar{v}}\right)\exp(-\phi),$$
 where
$\bar{v}=\frac{v_1v_2}{v_1+v_2}$, $\bar{a}=\frac{a_1v_1+a_2v_2}{v_1+v_2}$, 
$\phi=\frac{(a_1-a_2)^2}{2(v_1+v_2)}=\phi^{j+1}$.

We can also define a variable $\eta_{+}^{j+1}$ that is a
linear combination of $dX^{j+1}$, $dY^{j+1}$
and is
uncorrelated with $\eta^{j+1}$: 
$$\eta_{+}^{j+1}=\frac{-f_y\cdot dY^{j+1}+f_x\cdot dX^{j+1}}{\sqrt{f_x^2+f_y^2}}.$$
The observations do not affect $\eta_+^{j+1}$, so its mean and variance
are known. Given the means and variances of $\eta^{j+1}$, $\eta^{j+1}_{+}$ one can easily invert
the orthogonal matrix that connects them to  
$dX^{j+1}$, $dY^{j+1}$ and find the means and variances
$a_x,v_x$ of $dX^{j+1}$ and $a_y,v_y$ of
$dY^{j+1}$ after their modification by the observation 
(the subscripts on $a,v$ are labels, not differentiations).
Now one can produce values for $dX^{j+1},dY^{j+1}$:
$$dX^{j+1}=a_x+\sqrt{v_x}\xi_x,\quad  dY^{j+1}=a_y+\sqrt{v_y}\xi_y,$$
 where $\xi_x$, $\xi_y$
are the samples from $N(0,1)$ chosen at the beginning of the iteration. 
This completes the iteration.

This iteration converges to ${\bf X}^{n+1}$
such that $f({\bf X}^{n+1})=b^{n+1}+N(0,s)$,
and the phases $\phi^j$ converge to a limit $\phi=\phi_i$,
where the particle index $i$ has been restored. 
The time interval over which the solution is updated in each step is short,
and we do not expect any problem with convergence, either here or in the 
next section, and indeed there is none; in all cases the iteration converges
in a small number of steps. Note that after the iteration the variables $X^{n+1}_i,Y^{n+1}_i$ are
no longer independent- the observation creates a relation between them. 

Do this for all the particles. The particles are now samples of $P_{n+1}$, but they have been 
obtained by sampling different densities (remember that the parameters in the Gaussians
in equation (\ref{forward}) vary). One can get rid of this heterogeneity by viewing the factors
$\exp(-\phi)$ as weights and resampling,
 i.e., for each of $M$ random numbers
$\theta_k,k=1,\dots,M$ drawn from the uniform distribution on $[0,1]$,
choose a new ${\bf \hat X}^{n+1}_k={\bf X}^{n+1}_i$
such that
$Z^{-1}\sum_{j=1}^{i-1}\exp(-\phi_j) <  \theta_k \le Z^{-1}\sum_{j=1}^i \exp(-\phi_j)$
(where $Z=\sum_{j=1}^M \exp(-\phi_j)$), 
and then suppress the hat. 
We have traded the resampling of Bayesian filters for a resampling based on the 
normalizing factors of the several Gaussian densities; this is a worthwhile trade because
in a Bayesian filter one gets a set of samples many of which may have low
probability with respect to $P_{n+1}$, and here we have a set of samples each one of which
has high probability with respect to a pdf close to $P_{n+1}$. 

Note also that the resampling does not have to be done at every step- for example, one can
add up the phases for a given particle and resample only when the ratio of the
largest cumulative weight $\exp(-\sum \phi_i)$  to the smallest such weight exceeds some limit $L$ (the summation is over the  weights
accrued to a particular particle $i$ since the last resampling).
If one is worried by too many particles being close to 
each other ("depletion" in the 
Bayesian terminology), one can divide the set of 
particles into subsets of small size
and resample only inside those subsets, creating a 
greater diversity. As will be seen
in section 6, none of these strategies will be 
used here and we will resample fully at every step.  

\section{Backward sampling.}

The algorithm of the previous section is sufficient to create a  
filter, but accuracy may require an additional refinement. Every observation provides information
not only about the future but also about the past- it may, for example,
tag as improbable earlier states that had seemed probable before the observation was made; one
may have to go back and correct the past after every observation (this backward sampling is often misleadingly motivated solely by the need to create greater diversity
among the particles in a Bayesian filter).
As will be seen below, this backward sampling does not provide a significant boost to
accuracy in the present problem, but it is described here for the sake a completeness.

Given a set of particles at time $(n+1)\delta$, after a forward step and maybe a subsequent resampling,
one can figure out where each particle $i$  was in the previous two
steps, and have a partial history for each particle $i$:
${\bf X}_i^{n-1},{\bf X}_i^{n},{\bf X}_i^{n+1}$ (if resamples had
occurred,
some parts of that history may be shared among several current
particles). Knowing the first and the last member of this sequence, one can interpolate
for the middle term as in section 2, thus projecting information
backward. This requires that one recompute $d{\bf X}^n$. 

Let $d{\bf X}^{\text{tot}}=d{\bf X}^n+d{\bf X}^{n+1}$; in the present section this quantity is
assumed known and
remains fixed.
In the azimuth problem discussed here, one has to deal with the slight complication due to 
the fact that the mean of each increment is the value of the previous one, so that two successive
increments are related in a slightly more complicated way than usual. The displacement $dX^n$ is 
a $N(dX^{n-1},\sigma)$ variable, and $dX^{n+1}$ is a $N(dX^n,\sigma)$ variable, so that one goes from
$X^{n-1}$ to $X^{n+1}$ by sampling first a $(2dX^{n-1},4\sigma)$ variable that takes us 
from ${ X}^{n-1}$ to an intermediate point $P$, with a correction
by the observation half way up this first leg, and then one samples a 
$N(dX^{\text{tot}},\sigma)$ variable to reach $X^{n+1}$, and similarly for $Y$.
Let the variable that connects ${\bf X}^{n-1}$ to $P$ be
$d{\bf X}^{\text{new}}$,
so that what replaces $d{\bf X}^n$ is $d{\bf X}^{\text{new}}/2$.
Accordingly, we are looking for a new displacement
$d{\bf X}^{\text{new}}=(dX^{\text{new}},dY^{\text{new}})$, and for parameters
$a_x^{\text{new}},a_y^{\text{new}}, v_x^{\text{new}},v_y^{\text{new}}$ such that
\begin{align*}
&\exp\left(-\frac{\xi_x^2 +\xi_y^2}{2}\right)\nonumber\\
=&\exp\left(-\frac{(dX^{\text {new}}-2dX^{n-1})^2}{8\sigma}
-\frac{(dY^{\text {new}}-2dY^{n-1})^2}{8\sigma}\right)\\
&\times \exp\left(-\frac{(f^n-b^n)^2}{2s}\right)\nonumber\\
&\times
\exp\left(-\frac{(dX^{\text{new}}-dX^{\text{tot}})^2}{2\sigma}
-\frac{(dY^{\text{new}}-dX^{\text{tot}})^2}{2\sigma}\right)\exp(\phi)\nonumber\\
=&\exp\left(-\frac{(dX^{\text {new}}-\bar{a}_x)^2}{2v^{\text{new}}_x}
-\frac{(dY^{\text{new}}-\bar{a}_y)^2}{2v^{\text {new}}_y}\right), 
\end{align*}
where $f^n=f(X^{n-1}+dX^{\text{new}}/2,Y^{n-1}+dY^{\text{new}}/2)$ and 
$\xi_x$, $\xi_y$ are independent $N(0,1)$ Gaussian variables.
As in equation (\ref{forward}), the first equality embodies what we wish to 
accomplish-
find increments, functions of the reference variables, that sample the 
new pdf at time
$n\delta$ defined by the forward motion, the constraint imposed by the 
observation,
and by knowledge of the position at time $(n+1)\delta t$. The second 
equality states
that this is done by finding 
particle-dependent 
parameters for a Gaussian
density. 

We again find these parameters as well as the increments by iteration.
Much of the work is separate for the $X$ and $Y$ components
of the equations of  motion, so we write some of the equations for the 
$X$ component only.
Again set up an iteration for variables
$dX^{\text{new},j}=dX^j$ which converge to $dX^{\text{new}}$. Start with
$dX^0=0$. To find $dX^{j+1}$ given $dX^j$, approximate the observation 
equation (\ref{eq2}), as before, by
equation (\ref{obs}); define again variables $\eta^{j+1}, \eta^{j+1}_{+}$, one in the 
direction
of the approximate constraint and one orthogonal to it; in the direction
of the constraint multiply the pdfs as in the previous section; construct
new means $a^1_x,a^1_y$ and new variances
$v^1_x,v^1_y$ for $dX,dY$ at time $n$, taking into account the observation
at time $n$, again as before.
This also produces a phase $\phi=\phi_0$.

Now take into account that the location of the boat at time $n+1$ is
known;
this creates a new mean $\bar{a}_x$,
a new variance $\bar{v}_x$, and a new phase $\phi_x$, by
$\bar{v}=\frac{v_1v_2}{v_1+v_2}$, 
$\bar{a}_x=\frac{a_1v_1+a_2v_2}{v_1+v_2}$, 
$\phi_x=\frac{(a_1-a_2)^2}{v_1+v_2}$,
where $a_1=2a^1,v_1=4v^1_x,a_2=X^{\text{tot}},v_2=\sigma$.
Finally, find a new interpolated position $dX^{j+1}=a_x^{\text{new}}/2+\sqrt{v_x^{\text{new}}}\xi_x$ (the calculation for $dY^{j+1}$ is similar, with a 
phase $\phi_y$), and we are done.
The total phase for in this iteration is $\phi=\phi_0+\phi_x+\phi_y$.
 As the iterates $dX^j$ converge to $dX^{\text{new}}$, the phases converge
to a limit $\phi=\phi_i$. 
The probability of a particle arriving at the given position at time $(n+1)\delta t$ having been determined in the forward step,
there is no need to resample before comparing samples. 
Once one has the values of $\bf{X}^{\text{new}}$, a forward step gives corrected values of
${\bf X}^{n+1}$; one can use this interpolation process to correct estimates of ${\bf X}^k$ by  
subsequent observations 
for $k=n-1, k=n-2,\dots$, as many as are useful.

\section{Numerical results.}

Before presenting examples of numerical results for the azimuth problem, 
we  discuss
the accuracy one can expect. A single set of observations for our problem relies
on 160 samples of a $N(0,\sigma)$ variable. The maximum likelihood estimate of $\sigma$ given
these samples is a random variable with mean $\sigma$ and standard deviation $.11\sigma$.
We estimate the uncertainty in the position of the boat 
by picking a set of observations, then  making multiple
runs of the boat where the random components of the motion in the 
direction of the
constraint are frozen while the ones orthogonal to it are sampled 
over and over from
the suitable Gaussian density, then computing the distances to the fixed 
observations, estimating the
standard deviation of these differences, and accepting the trajectory 
if the estimated standard deviation is  
within one standard deviation of the nominal value of $s$. This process generates a family of boat 
trajectories compatible with the given observations. In Table I 
we display the standard
deviations of the differences between the resulting paths and the 
original path that produced
the observations after the number of steps indicated there 
(the means of these differences
are statistically indistinguishable from zero). 
This Table provides an estimate of
the accuracy we can expect. 
It is fair to assume that these standard
deviations are underestimates of the uncertainty- a variation of a single standard deviation in $s$ is a 
strict constraint, and we allowed no variation in $\sigma$.

\begin{center}
Table I\\
Intrinsic uncertainty in the azimuth problem
\end{center}

\begin{center}
\begin{tabular}{|c|c|c|}
\hline
step&            $x$ component&                     $y$  component\\
\hline
40&              .0005                                 & .21 \\                      
80&              .004                                  &  .58 \\                      
120&             .010                                  &  .88 \\
160&             .017                                  &  .95 \\
\hline
\end{tabular}
\end{center}

\vspace{12pt}
If one wants reliable information about the performance of the filter, it is not sufficient to run the boat once, record observations,
and then use the filter to reconstruct the boat's path, because the difference
between the true path and the reconstruction is a random variable which may be
accidentally atypically small or atypically large. We have therefore run a large number of such  
reconstructions and computed the means and standard deviations of the discrepancies 
between path and reconstruction as a function of the number of steps and of other
parameters. 
In Tables II and III we display the means and standard deviations of these discrepancies (not of their mean!) in the 
the x and y components of the paths with 2000 runs, at the steps and numbers of particles
indicated, with no backward sampling. 
(Ref. \cite{Gi1} used  100 particles).
On the average the error is zero, and the error that can be expected
in any one run is of the order of magnitude of the unavoidable error. 
The standard 
deviation of the discrepancy is not significantly smaller with 2 particles that with 100-
the main source of the discrepancy is the uncertainty in the data.
Most of time one single particle (no resampling) is enough; however,
a single particle may temporarily stray into low-probability areas and creates large
arguments and numerical difficulties in the various functions used in the program. 
Two particles with resampling keep each other within bounds, because if one of them 
strays it gets replaced by a
replica of the other.
The various more sophisticated resampling strategies at the end of section 4 make no discernible difference here, and backward
sampling does not help much either,
because they too are unable to remedy the limitation of the data set.


\begin{center}
     Table IIa\\

Mean and standard variation of the discrepancy between synthetic data and
their reconstruction, 2000 runs, no back step, 100 particles
\end{center}
\begin{center}
\begin{tabular}{|c|c|c|c|c|}\hline
n. of steps &    \multicolumn{2}{c|}{x component} & \multicolumn{2}{c|}{ y component}\\\hline

          &   mean  &  s.d. &              mean & s.d.\\\hline

40       &    .0004 &  .04 &             .0001  & .17  \\\hline

80        &  -.001  &  .04  &           -.01  &   .54\\\hline

120       &  -.0008 &  .07   &          -.03  &  1.02\\\hline

160     &    -.002  &  .18    &         -.05  &  1.56 \\\hline
\end{tabular}
\end{center}

\begin{center}
     Table IIb\\

Mean and standard variation of the discrepancy between synthetic data and
their reconstruction, 2000 runs, no back step, 2  particles
\end{center}
\begin{center}
\begin{tabular}{|c|c|c|c|c|}\hline
n. of steps &    \multicolumn{2}{c|}{x component} & \multicolumn{2}{c|}{ y component}\\\hline

          &   mean  &  s.d. &              mean & s.d.\\\hline
                                                                                                
40      &   .002 &     .17 &            -.0004 &   .20\\\hline

80      &   .01   &    .43  &           -.0006  &  .58\\\hline

120     &   .01    &   .57   &           .009    & 1.08\\\hline

160      &  .006    &  .54    &          .01     & 1.67\\\hline
\end{tabular}
\end{center}

In Figure 1 we plot a sample boat path, its reconstruction, and the reconstructions obtained 
(i) when the initial data for the reconstruction are strongly perturbed
(here, the initial data
for $x,y$ were perturbed initially by, respectively, $.1$ and $.4$), and (ii) when the
value of $\sigma$ assumed in the reconstruction is random:
$\sigma=N(\sigma_0,\epsilon\sigma_0)$, where $\sigma_0$ is the constant value used until now and
$\epsilon=0.4$ but the calculation is otherwise identical. This produces variations in $\sigma$ of  the order of
$40\%$; any larger variance in the perturbations produced negative value of $\sigma$. 
The differences between the reconstructions and the true path remain within
the acceptable range of errors. These graphs show that the filter has little sensitivity to 
perturbations (we did not calculate statistics here because the insensitivity
holds for each individual run).

\begin{figure}

\begin{center}

\includegraphics[scale=0.3]{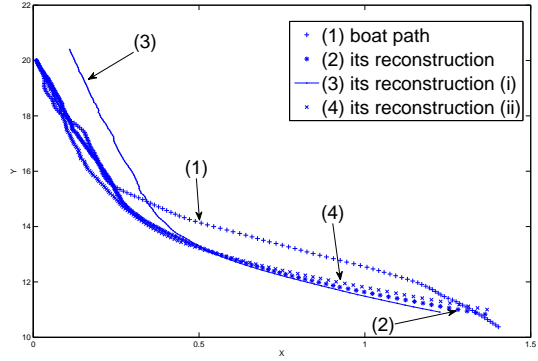}
\caption{Some boat trajectories (explained in the text)}




\end{center}

\end{figure}

We now estimate the parameter $\sigma$ from data. The filter needs an estimate
of $\sigma$ to function, 
call this estimate $\sigma_{\text{assumed}}$. If $\sigma_{\text{assumed}}\ne \sigma$,
the other assumptions used to produce the data set (e.g. independence of the displacements and of the observations) are also false, and all one has to do
is detect the fallacy. We do it by picking a trajectory of a particle and computing
the quantity
$$D=
\frac{(\sum_2^J(dX^{j+1}-dX^j))^2+(\sum_2^J(dY^{j+1}-dY^j))^2}{
\sum_2^J(dX^{j+1}-dX^j)^2+\sum_2^J(dY^{j+1}-dY^j)^2}.
$$
If the increments are independent 
 then on the average $D=1$; we will try to find the real $\sigma$ by finding a value of $\sigma_{\text{assumed}}$
for which this happens. We chose $J=40$ (the early part of a trajectory is less noisy than
the later parts).

As we already know, a single run cannot provide an accurate estimate of $\sigma$, and accuracy
in the reconstruction depends on how many runs are used. In Table III we display some
values of $D$ averaged over 200 and over 5000 runs as a function of the ratio of
 $\sigma_{\text{assumed}}$ to the value of $\sigma$ used to generate the data. From the longer computation one can find the correct value of $\sigma$  
with an error of about $3\%$, while with 200 runs the uncertainty is about $10\%$.

\begin{center}
    Table III\\

The mean of the discriminant D as a function of $\sigma_{\text{assumed}}/\sigma$, 30 particles 
\end{center}

\begin{center}
\begin{tabular}{|c|c|c|}\hline

$ \sigma_{\text{assumed}}/\sigma$  &  5000 runs  &        200 runs   \\\hline
.5              &   1.14 $\pm$ .01  &     1.21 $\pm$ .08\\\hline

.6             &    1.08 $\pm$ .01   &    1.14 $\pm$ .07    \\\hline              

.7              &   1.05 $\pm$ .01   &    1.10 $\pm$ .07\\\hline

.8               &  1.04 $\pm$ .01   &     1.14 $\pm$ .07\\\hline

.9                & 1.00 $\pm$ .01   &     1.01 $\pm$ .07\\\hline

1.0              &  1.00 $\pm$ .01   &      .96 $\pm$ .07\\\hline

1.1              &   .97 $\pm$ .01   &      1.01 $\pm$ .07\\\hline

1.2              &   .94 $\pm$ .01   &       .99 $\pm$ .07\\\hline

1.3              &   .93 $\pm$  .01   &     1.02 $\pm$ .07   \\\hline                    

1.4              &   .90 $\pm$ .01     &     .85 $\pm$ .06  \\\hline              

1.5              &   .89  $\pm$ .01     &    .93 $\pm$ .07\\\hline

2.0              &   .86  $\pm$ .01     &    .78 $\pm$ .05\\\hline
\end{tabular}
\end{center}

\section{Conclusions.}

We have exhibited a non-Bayesian filtering method, related to recent work on chainless sampling, designed to focus particle paths more sharply and thus require fewer of them, at the cost of
an added complexity in the evaluation of each path. The main features of the algorithm are
a representation of a new pdf by means of a set of functions of Gaussian variables  and a resampling based on normalization
factors. The construction was demonstrated on a standard ill-conditioned test problem.
Further applications will be published elsewhere.

\section{Acknowledgments.}
We would like to thank Prof. R. Kupferman, Prof. R. Miller, and  Dr. J. Weare
for asking searching questions and providing good advice. 
This work was supported in part by the Director, Office of Science,
Computational and Technology Research, U.S.\ Department of Energy under
Contract No.\ DE-AC02-05CH11231, and by the National Science Foundation under grant
 DMS-0705910.

\end{document}